\documentclass[12pt]{article}
\usepackage{amssymb}
\usepackage{amsfonts}
\usepackage{amsmath}
\usepackage{amsthm}
\usepackage{amscd}
\usepackage{eucal}
\usepackage{latexsym}
\usepackage{fancyhdr}
\theoremstyle{plain}
\newtheorem{thm}{Theorem}

\newtheorem{conj}{Conjecture}
\theoremstyle{remark}
\newtheorem*{remark}{Remark}
\newtheorem*{remarks}{Remarks}

\numberwithin{thm}{section}
\numberwithin{prop}{section}
\begin{document}
\title{Quasi-symmetric functions and mod $p$ multiple harmonic sums}
\author{Michael E. Hoffman\\
\small Dept. of Mathematics\\[-0.8ex]
\small U. S. Naval Academy, Annapolis, MD 21402\\[-0.8ex]
\small \texttt{meh@usna.edu}}
\pagestyle{fancy}
\lhead{M. E. Hoffman}
\rhead{QSym and mod $p$ multiple harmonic sums}
\date{\small \today
}
\maketitle
\def\al{\alpha}
\def\be{\beta}
\def\de{\delta}
\def\ep{\epsilon}
\def\zt{\zeta}
\def\zts{\zeta^{\star}}
\def\la{\lambda}
\def\si{\sigma}
\def\Ga{\Gamma}
\def\De{\Delta}
\def\Si{\Sigma}
\def\del{\nabla}
\def\tilde{\widetilde}
\def\<{\langle}
\def\>{\rangle}
\def\Qxy{\mathbf Q\< x,y\>}
\def\Zxy{\mathbf Z\< x,y\>}
\def\RN{\mathbf R^{\mathbf N}}
\def\Sym{\operatorname{Sym}}
\def\QSym{\operatorname{QSym}}
\def\id{\operatorname{id}}
\def\card{\operatorname{card}}
\def\h{\operatorname{ht}}
\def\Zp{\mathbf Z/p\mathbf Z}
\def\A{\mathcal A}
\def\Op{\mathcal O}
\def\H{\mathfrak H}
\def\P{\mathfrak P}
\def\stir#1#2{\genfrac{\{}{\}}{0pt}{}{#1}{#2}}
\def\a{\zts_{5211}}
\def\b{\zts_{6111}}
\def\c{\zts_{81}}
\def\e{\zts_{21}\zts_{411}}
\newcommand{\stone}[2]{\genfrac{[}{]}{0pt}{}{#1}{#2}}
\newcommand{\sttwo}[2]{\genfrac{\{}{\}}{0pt}{}{#1}{#2}}
\begin{abstract}
We present a number of results about (finite) multiple harmonic sums
modulo a prime, which provide interesting parallels to known
results about multiple zeta values (i.e., infinite multiple
harmonic series).
In particular, we prove a ``duality'' result
for mod $p$ harmonic sums similar to (but distinct from) that for 
multiple zeta values.
We also exploit the Hopf algebra structure of the quasi-symmetric
functions to do calculations with multiple harmonic sums mod $p$,
and obtain, for each weight through 9, a set of generators for the
space of weight-$n$ multiple harmonic sums mod $p$.
When combined with recent work, the results of this paper offer
significant evidence that the number of quantities needed to 
generate the weight-$n$ multiple harmonic sums mod $p$ is the
$n$th Padovan number (OEIS sequence A000931).
\end{abstract}
\section{Introduction}
\par
In recent years there has been considerable interest in the
multiple zeta values
\begin{equation}
\label{mzv}
\zt(i_1,i_2,\dots,i_k)=\sum_{n_1>n_2>\dots>n_k\ge 1} \frac1{n_1^{i_1}
n_2^{i_2}\cdots n_k^{i_k}} ,
\end{equation}
which converge for $i_1>1$; see, e.g., \cite{BB,H1,H2,HO,H4,W,Z}.  The
multiple zeta values are limits of the finite harmonic sums
\begin{equation}
\label{fmzv}
\zt_n(i_1,i_2,\dots,i_k)=\sum_{n\ge n_1>n_2>\dots>n_k\ge 1} \frac1{n_1^{i_1}
n_2^{i_2}\cdots n_k^{i_k}} ,
\end{equation}
as $n\to\infty$.  In fact, these finite sums themselves are of interest 
in physics \cite{BK,Ve}.  Both types of quantities are indexed by the
exponent string $(i_1,\dots,i_k)$:  we refer to $k$ as the length and
$i_1+\dots+i_k$ as the weight.
\par
The multiple zeta values (\ref{mzv}) form an algebra under multiplication.  
In fact, 
there are two distinct multiplications on the set of multiple
zeta values, the harmonic (or ``stuffle'') product and the shuffle product, 
coming from the representation of multiple zeta values as iterated series
and as iterated integrals respectively (see \S4 of \cite{HO} for a discussion).
Only the first of these products applies to the finite sums (\ref{fmzv}).
\par
In this paper we develop a theory of the mod $p$ values ($p$ a prime),
of the finite harmonic sums (2) when $n=p-1$.  Our main tools are the 
structure of the algebra $\QSym$ of quasi-symmetric functions (which 
formalizes the harmonic product) and a mod $p$ duality result (Theorem 
4.7 below) which corresponds to the duality theorem for multiple zeta 
values (see, e.g., Corollary 6.2 of \cite{H2}).  
A brief account of the results of \S4 of this paper appeared earlier 
in the last section of \cite{H4}.
\par
The theory of mod $p$ finite harmonic sums has several interesting
differences and similarities to the theory of multiple zeta values.
While the values of (1) with length $k=1$ are the classical zeta values
$\zt(i_1)$, it is well known that the mod $p$ value of 
$$
\zt_{p-1}(i_1)=\sum_{j=1}^{p-1}\frac1{j^{i_1}}
$$
is zero when $p>i_1+1$.  But for $k=2$ we have the relation
$$
\zt_{p-1}(i_1,1)\equiv B_{p-i_1-1} \mod p ,
$$
with Bernoulli numbers (see Theorem 6.1 below).  In the case of 
multiple zeta values,
it appears that multiple zeta values not expressible in terms of 
the classical zeta values occur first in weight 8 (e.g., $\zt(6,2)$);
for the finite harmonic sums mod $p$, the first ones not expressible 
in terms of Bernoulli numbers appear to occur in weight 8 (i.e., 
$\zt_{p-1}(6,1,1)$).
\par
As we have already mentioned, there is a duality result in both
theories.  If we code the exponent strings $(i_1,i_2,\dots,i_k)$
occurring in (1) and (2) by the monomial $x^{i_1-1}y\cdots x^{i_k-1}y$
in noncommuting variables $x$ and $y$, then the duality map for
multiple zeta values is the anti-automorphism of $\Qxy$ that
exchanges $x$ and $y$, while that for the mod $p$ finite harmonic sums
is the automorphism $\psi$ of $\Zxy$ with $\psi(x)=x+y$ and $\psi(y)=-y$.
\par
The dimension of the rational vector space of multiple zeta
values of weight $n$ was conjectured by D. Zagier \cite{Z} to be the
numbers $d_n$ given by $d_1=0$, $d_2=d_3=1$, and $d_n=d_{n-2}+d_{n-3}$.
It has been proved by T. Terasoma \cite{T} and A. B. Goncharov \cite{G} 
that the true dimension of this rational vector space is at most $d_n$.
The corresponding numbers for the finite harmonic sums (2) are
the minimal numbers $c_n$ of weight-$n$ harmonic sums needed to
generate all weight-$n$ harmonic sums mod $p$ for $p>n+1$.
From the calculations of \S7 below, supplemented by the recent
results of Kh. Hessami Pilehrood, T. Hessami Pilehrood and 
R. Tauraso \cite{HPT}, the first few values of $c_n$ are the 
following.
\par
\begin{center}
\begin{tabular}{|c|c|c|c|c|c|c|c|c|c|}\hline
$n$ & 1 & 2 & 3 & 4 & 5 & 6 & 7 & 8 & 9 \\ \hline
$c_n$ & 0 & 0 & 1 & 0 & 1 & 1 & 1 & 2 & 2\\ \hline
\end{tabular}
\end{center}
We make the following conjecture.
\begin{conj}
The numbers $c_n$ are the Padovan numbers (OEIS sequence A000931),
i.e. the sequence with $c_1=c_2=0$, $c_3=1$ and $c_n=c_{n-2}+c_{n-3}$.
\end{conj}
If this conjecture is true, it follows that $c_n=d_{n-3}$ for $n\ge 4$,
as follows from consideration of the generating functions.
\par
This paper is organized as follows.  In \S2 we give an exposition
of the algebra $\QSym$ of quasi-symmetric functions.  In \S3 we
describe the Hopf algebra structure of $\QSym$, particularly
the antipode, and introduce some integral bases of $\QSym$
important in the sequel.  In \S4 we define bases $\zt_n(I)$ and $\zts_n(I)$
for the set of finite harmonic sums, and prove some basic
results about the mod $p$ case.  In \S5 we discuss a particular
class of multiple harmonic sums, the height one sums, and express
them in terms of Stirling numbers of the first and second kinds.
In \S6 we obtain some results about $\zts_{p-1}(I)$ mod $p$ when $I$ 
has short length.
Finally, in \S7 we apply the preceding results to find sets of 
products of height one sums that generate all multiple harmonic
sums mod $p$ through weight 9.
\par
Except for the introduction and some references to later work, this
paper is the one written during a stay at the Max Planck Institut f\"ur 
Mathematik (MPIM) in Bonn in 2003-04, and disseminated as the preprint
MPI 04-03.
The notation has been modified to bring it in line with more recent 
papers, and some minor mistakes in the preprint have been corrected.
While many of the results obtained here appeared in the work of 
J. Zhao \cite{Zh}, the author believes the approach and point of view
herein differs significantly from other published works.
\par
The author thanks the directors of the MPIM for their hospitality,
and also thanks Joe Buhler for his prompt response to a query about 
irregular primes.
\section{Harmonic Algebra and Quasi-Symmetric Functions}
\par
Let $\H$ be the underlying graded abelian group of the noncommutative
polynomial algebra $\Zxy$, where $x$ and $y$ both have degree 1.
For a word $w$ of $\H$, we refer to the total degree of $w$ as
its weight (denoted $|w|$) and its $y$-degree as its length
(denoted $\ell(w)$).
We can define a commutative multiplication $*$ on $\H$ by requiring
that it distribute over the addition and that it satisfy the following
axioms:
\begin{itemize}
\item[H1.]
For any word $w$, $1*w=w*1=w$;
\item[H2.]
For any word $w$ and integer $n\ge 1$, $x^n*w=w*x^n=wx^n$;
\item[H3.]
For any words $w_1,w_2$ and integers $p,q\ge 0$,
$$
x^pyw_1*x^qyw_2=x^py(w_1*x^qyw_2)+x^qy(x^pyw_1*w_2)+x^{p+q+1}y(w_1*w_2) .
$$
\end{itemize}
Note that axiom (H3) allows the $*$-product of any pair of words
to be computed recursively, since each $*$-product on the right has
fewer factors of $y$ than the $*$-product on the left-hand side.
Induction on length establishes the following.
\begin{thm} The $*$-product is commutative and associative.
\end{thm}
We refer to $\H$ together with its commutative multiplication $*$
as the (integral) harmonic algebra $(\H,*)$.  
Let $\H^1$ be the additive subgroup $\mathbf Z1 +\H y$ of $\H$; it
is evidently a subalgebra of $(\H,*)$.  
Note that any word $w\in\H^1$ can be written in terms of
the elements $z_i=x^{i-1}y$, and that the length $\ell(w)$
is the number of factors $z_i$ of $w$ when expressed this way.
We can rewrite the inductive rule (H3) for the $*$-product as
\begin{equation}
z_pw_1*z_qw_2=z_p(w_1*z_qw_2)+z_q(z_pw_1*w_2)+z_{p+q}(w_1*w_2).
\label{zrule}
\end{equation}
\par
Now for each positive integer $n$, define a homomorphism $\phi_n:\H^1\to
\mathbf Z[t_1,\dots,t_n]$ of graded abelian groups (where $|t_i|=1$ for all $i$) 
as follows.  Let $\phi_n(1)=1$ and 
$$
\phi_n(z_{i_1}z_{i_2}\cdots z_{i_k})=
\sum_{1\le n_1<n_2<\dots<n_k\le n}t_{n_1}^{i_1}t_{n_2}^{i_2}
\cdots t_{n_k}^{i_k}
$$
for words of length $k\le n$, and let $\phi(w)=0$ for words of
length greater than $n$; extend $\phi_n$ additively to $\H^1$.
Because the rule (\ref{zrule}) corresponds to multiplication of series,
$\phi_n$ is an algebra homomorphism of $(\H^1,*)$ into $\mathbf Z[t_1,\dots,t_n]$, 
and $\phi_n$ is evidently injective through degree $n$.  For each $m\ge n$, there 
is a restriction map
$$
\rho_{m,n}:\mathbf Z[t_1,\dots,t_m]\to \mathbf Z[t_1,\dots,t_n]
$$
defined by
$$
\rho_{m,n}(t_i)=\begin{cases} t_i,& i\le n\\ 0,& i>n . \end{cases}
$$
The inverse limit
$$
\P=\projlim_n \mathbf Z[t_1,\dots,t_n]
$$
is the subalgebra of $\mathbf Z[[t_1,t_2,\dots]]$ consisting of
those formal power series of bounded degree.  The $\phi_n$
commute with the restriction maps, so they define a homomorphism
$\phi:\H^1\to\P$.
\par
Inside $\P$ is the algebra of symmetric functions
$$
\Sym =\projlim_n \mathbf Z[t_1,\dots,t_n]^{\Si_n}
$$
and also the algebra $\QSym$ of quasi-symmetric functions (first
described in \cite{Ges}).  We say
a formal series $p\in\P$ is in $\QSym$ if the coefficient of 
$t_{i_1}^{p_1}\cdots t_{i_k}^{p_k}$ in $p$ is the same as the
coefficient of $t_{j_1}^{p_1}\cdots t_{j_k}^{p_k}$ in $p$ 
whenever $i_1<i_2<\dots<i_k$ and $j_1<j_2<\dots<j_k$.  Evidently
$\Sym\subset\QSym$.  An integral basis for $\QSym$ is given by
the monomial quasi-symmetric functions 
$$
M_{(p_1,p_2,\dots,p_k)}=\sum_{i_1<i_2<\dots<i_k}t_{i_1}^{p_1}t_{i_2}^{p_2}
\cdots t_{i_k}^{p_k} ,
$$
which are indexed by compositions $(p_1,\dots,p_k)$.  Since 
$\phi(z_{i_1}\cdots z_{i_k})=M_{(i_1,\dots,i_k)}$, we have the following
result.
\begin{thm} $\phi$ is an isomorphism of $(\H^1,*)$ onto $\QSym$.
\end{thm}
\par
As is well known, the algebra $\Sym$ of symmetric functions is generated
by the elementary symmetric functions $e_i$, as well as by the complete
symmetric functions $h_i$.  The power-sum symmetric functions $p_i$
generate $\Sym\otimes\mathbf Q$, but only generate a subalgebra of
$\Sym$ over $\mathbf Z$.  
It is easy to see that $\phi(z_1^i)=e_i$ and $\phi(z_i)=p_i$.
\par
A monomial quasi-symmetric function $M_I$ is in $\Sym$ exactly when
all parts of $I$ are the same.  Given $I=(i_1,\dots,i_k)$, the symmetric
group $\Si_k$ on $k$ letters acts on $I$ by $\si\cdot I=(i_{\si(1)},\dots,
i_{\si(k)})$, and the symmetrization
\begin{equation}
\sum_{\si\in\Si_k} M_{\si\cdot I}
\label{sym}
\end{equation}
of $M_I$ is evidently in $\Sym$.  Hence the element (\ref{sym}) can be 
written as a sum of rational multiples of the power-sums $p_i=M_{(i)}$.
But in fact (\ref{sym}) is a sum of products of the $M_{(i)}$ with
integral coefficients that can be given explicitly as follows.
\begin{thm}
Let $I=(i_1,i_2,\dots,i_k)$ be a composition.  Then
\begin{equation*}
\sum_{\si\in\Si_k} M_{\si\cdot I}=\sum_{\text{partitions 
$\mathcal B=\{B_1,\dots,B_l\}$ of 
$\{1,\dots,k\}$}}(-1)^{k-l}c(\mathcal B)M_{(b_1)}M_{(b_2)}\cdots M_{(b_l)} ,
\end{equation*}
where $c(\mathcal B)=(\card B_1-1)!\cdots (\card B_l-1)!$ and
$b_s=\sum_{j\in B_s} i_j$.
\end{thm}
\begin{proof} 
While the argument used to prove Theorem 2.2 of \cite{H1} can be adapted
to prove this result, we shall use the M\"obius inversion formula as follows.  
The set $\Pi_k$ of partitions of the set $\{1,\dots,k\}$,
with the partial order given by refinement, is a finite semimodular
lattice.  We also have, for any $\mathcal C=\{C_1,\dots,C_p\}\in\Pi_k$,
$$
M_{(c_1)}M_{(c_2)}\cdots M_{(c_p)}=\sum_{\mathcal B=\{B_1,\dots,B_l\}\preceq
\mathcal C}\sum_{\si\in\Si_l}M_{\si\cdot (b_1,\dots,b_l)},
$$
where $b_s=\sum_{j\in B_s}i_j$ and $c_s=\sum_{j\in C_s}i_j$.  Then
the M\"obius inversion formula gives
\begin{equation}
\sum_{\si\in\Si_p}M_{\si\cdot (c_1,\dots,c_p)}=\sum_{\mathcal 
B=\{B_1,\dots,B_l\}\preceq\mathcal C} \mu(\mathcal B,\mathcal C)
M_{(b_1)}\cdots M_{(b_l)} .
\label{mif}
\end{equation}
The conclusion follows by taking $\mathcal C=\{\{1\},\{2\},\dots,\{k\}\}$
in equation (\ref{mif}) and noting that in this case 
$\mu(\mathcal B,\mathcal C)$=$(-1)^{k-l}c(\mathcal B)$ 
(see \cite{St}, Example 3.10.4).
\end{proof}
\par
The rational algebra of quasi-symmetric functions $\QSym\otimes\mathbf Q$
was shown to be a polynomial algebra by C. Malvenuto and C. Reutenauer 
\cite{MR}.
To describe the generators, we put an order on the words of $\H$ by
setting $x<y$ and extending it lexicographically.  A word is called
Lyndon if it is smaller than any of its proper right factors, i.e.
for $w$ Lyndon we have $w<v$ whenever $w=uv$ for $u,v\ne 1$.  
Let $\mathcal L$ be the set of Lyndon words in $\H^1$.
Then the following result was proved in \cite{MR}.
\begin{thm}
$\QSym\otimes\mathbf Q$ is a polynomial algebra on the generators
$\phi(w)$, where $w\in \mathcal L$.
\end{thm}
The integral structure is more subtle.  The result is again that
the algebra is polynomial, but one has to use a different set
of generators.  
Let $\mathcal L^{mod}$ be the set of modified Lyndon words in $\H^1$,
i.e., whenever $w=z_{i_1}^{p_1}\cdots z_{i_k}^{p_k}\in\mathcal L$
has a common factor $d$ dividing all the exponents $p_i$, replace
it with $v^d$, where $v=z_{i_1}^{p_1/d}\cdots z_{i_k}^{p_k/d}$.
The following result was first stated by E. J. Ditters, though the first 
correct proof seems to be due to M. Hazewinkel \cite{Hz}.  
\begin{thm}
$\QSym$ is a polynomial algebra on the generators $\phi(w)$, where
$w\in\mathcal L^{mod}$.
\end{thm}
\section{Quasi-Symmetric Functions as a Hopf Algebra}
\par
For definitions and basic results on Hopf algebras we refer the reader
to \cite{MM}, \cite{S}, and \cite{Ka}.
The algebra $(\H^1,*)\cong\QSym$ has a Hopf algebra structure with coproduct
$\De$ defined by
\begin{equation}
\De(z_{i_1}z_{i_2}\cdots z_{i_n})=
\sum_{j=0}^n z_{i_1}\cdots z_{i_j}\otimes z_{i_{j+1}}\cdots z_{i_n} ,
\label{comul}
\end{equation}
and counit $\ep$ with $\ep(u)=0$ for all elements $u$ of positive degree.
This extends the well-known Hopf algebra structure on the algebra
$\Sym$ (as described in \cite{Gei}), in which the elementary symmetric
functions $e_i$ ($\leftrightarrow y^i$) and complete symmetric functions
$h_i$ are divided powers, while the power sums $p_i$ ($\leftrightarrow z_i$)
are primitive.  The Hopf algebra $(\H^1,*,\De)$ is commutative but not
cocommutative.  Its (graded) dual is the Hopf algebra of noncommutative
symmetric functions as defined in \cite{Gel}.
\par
Now $\QSym$ has various integral bases besides the $M_I$, also
indexed by compositions.  For compositions $I,J$,
we say $I$ refines $J$ (denoted $I\succ J$) if $J$ can be obtained
from $I$ by combining some of its adjacent parts.  Then
the fundamental quasi-symmetric functions are given by
\begin{equation}
F_I=\sum_{J\succeq I} M_J,
\label{fun}
\end{equation}
and the essential quasi-symmetric functions are given by
\begin{equation}
E_I=\sum_{J\preceq I} M_J .
\label{ess}
\end{equation}
\par
It will be useful to have some additional notations for compositions.
We adapt the notation used the previous section for words, so 
for $I=(i_1,\dots,i_k)$ the weight of $I$ is $|I|=i_1+\dots+i_k$,
and $k=\ell(I)$ is the length of $I$.  
For $I=(i_1,\dots,i_k)$, the reversed composition $(i_k,\dots,i_1)$ will
be denoted $\bar I$:  of course reversal preserves weight, length and
refinement (i.e., $I\succeq J$ implies $\bar I\succeq\bar J$).  We
write $I\sqcup J$ for the juxtaposition of $I$ and $J$, so $\overline
{I\sqcup J}=\bar J\sqcup\bar I$.
\par
Compositions of weight $n$ are in 1-to-1 correspondence with subsets of 
$\{1,2,\dots,n-1\}$ via partial sums
\begin{equation}
(i_1,i_2,\dots,i_k)\to \{i_1,i_1+i_2,\dots,i_1+\dots+i_{k-1}\},
\label{corr}
\end{equation}
and this correspondence is an isomorphism of posets, i.e.,
$I\succeq J$ if and only if the subset corresponding to $I$ contains
that corresponding to $J$.  Complementation in the power set thus
gives rise to an involution $I\to I^*$; e.g., $(1,1,2)^*=(3,1)$.
Evidently $|I^*|=|I|$ and $\ell(I)+\ell(I^*)=|I|+1$.  Also, $I\preceq J$
if and only if $I^*\succeq J^*$.  The complementation operation
commutes with reversal, so the notation $\bar I^*$ is unambiguous.
\par
Because of the correspondence (\ref{corr}) between the poset of 
compositions
of $n$ and the poset of subsets of $\{1,2,\dots,n-1\}$, it follows
that the M\"obius function for compositions of $n$ is given by
$$
\mu(I,J)=(-1)^{\ell(I)-\ell(J)} .
$$
Thus, e.g., the M\"obius inversion formula applied to equation (\ref{ess})
is
\begin{equation}
M_I=\sum_{J\preceq I}(-1)^{\ell(I)-\ell(J)}E_J .
\label{sse}
\end{equation}
\par
Since $\QSym$ is a commutative Hopf algebra, its antipode
$S$ is an automorphism of $\QSym$ and $S^2=\id$.  Now $S$ can
be given by the following explicit formulas:  for proof see
\cite{Eh} or \cite{H3}.
\begin{thm} The antipode $S$ of $\QSym$ is given by 
\begin{itemize}
\item[1.]
$S(M_I)=\sum_{I_1\sqcup I_2\sqcup\cdots\sqcup I_l=I} (-1)^l M_{I_1}M_{I_2}
\cdots M_{I_l};$
\item[2.]
$S(M_I)=(-1)^{\ell(I)}E_{\bar I}$.
\end{itemize}
\end{thm}
Part (2) of this result implies a number of facts about the $E_I$.
First, the $E_I$ have almost the same 
multiplication rules as the $M_I$:  if $T$ is the automorphism of 
$\QSym$ sending $M_I$ to $M_{\bar I}$, then $S T$ takes  any identity 
among the $M_I$ to an identity among the $E_I$ that differs only in signs.  
For example, since
$$
M_{(2)}M_{(3)}=M_{(2,3)}+M_{(3,2)}+M_{(5)}
$$
we have
$$
E_{(2)}E_{(3)}=E_{(2,3)}+E_{(3,2)}-E_{(5)} .
$$
Second, the coproduct formula (\ref{comul}), which in terms of the $M_I$ says
$$
\De(M_I)=\sum_{J\sqcup K=I} M_J\otimes M_K ,
$$
can also be written
$$
\De(E_I)=\sum_{J\sqcup K=I} E_J\otimes E_K,
$$
because of the standard Hopf algebra relation $\De S=(S\otimes S)\De^{op}$.
Finally, if we apply $S$ to both sides of Theorem 2.3 we obtain,
for any composition $I=(i_1,\dots,i_k)$,
\begin{equation}
\sum_{\si\in\Si_k} E_{\si\cdot I}=\sum_{\text{partitions 
$\mathcal B=\{B_1,\dots,B_l\}$ of 
$\{1,\dots,k\}$}}c(\mathcal B)M_{(b_1)}M_{(b_2)}\cdots M_{(b_l)} ,
\label{EM}
\end{equation}
where as above $c(\mathcal B)=(\card B_1-1)!\cdots (\card B_l-1)!$ and
$b_s=\sum_{j\in B_s}i_j$.
\par
Now define an automorphism $\psi$ of $\Zxy$ by 
\begin{equation}
\psi(x)=x+y,\quad \psi(y)=-y .
\label{psi}
\end{equation}
Evidently $\psi^2=\id$, and $\psi(\H^1)=\H^1$.
Thus $\psi$ defines an additive involution of $\H^1\cong
\QSym$ (which is {\it not}, however, a homomorphism for the 
$*$-product).  We can describe the action of $\psi$ on the integral
bases for $\QSym$ as follows.
\begin{thm} For any composition $I$,
\begin{itemize}
\item [1.]
$\psi(M_I)=(-1)^{\ell(I)}F_I$;
\item [2.]
$\psi(E_I)=-E_{I^*}$.
\end{itemize}
\end{thm}
\begin{proof} Suppose $w=w(I)$ is the word in $x$ and $y$ corresponding
to a composition $I$.  Then evidently substituting $y$ in place of any
particular factor $x$ in $w$ corresponds to splitting a part of $I$.  
With this observation, part (1) is clear (there is also one factor of 
$-1$ for each occurrence of $y$ in $w$).
\par
Now we prove part (2).  We have
$$
\psi(E_I)=\sum_{J\preceq I}\psi(M_J)=
\sum_{J\preceq I} (-1)^{\ell(J)}F_J
$$
from part (1).  From Example 1 of \cite{H3}, $S(F_I)=(-1)^{|I|}
F_{\bar I^*}$, where $S$ is the antipode of $\QSym$.  Thus
\begin{multline*}
S\psi(E_I)=\sum_{J\preceq I}(-1)^{\ell(J)+|J|}F_{\bar J^*}
=-\sum_{J\preceq I}(-1)^{\ell(J^*)}F_{\bar J^*}\\
=-\sum_{\bar J^*\succeq\bar I^*}(-1)^{\ell(\bar J^*)}F_{\bar J^*}
=-\sum_{K\succeq\bar I^*}(-1)^{\ell(K)}F_K .
\end{multline*}
Now by M\"obius inversion of equation (\ref{fun})
$$
M_I=\sum_{I\preceq J}(-1)^{\ell(I)-\ell(J)}F_J ,
$$
and so
$$
S\psi(E_I)=-(-1)^{\ell(\bar I^*)} M_{\bar I^*}
$$
Apply $S$ be both sides to get
$$
\psi(E_I)=-(-1)^{\ell(I^*)}(-1)^{\ell(\bar I^*)}E_{I^*}=-E_{I^*}.
$$
\end{proof}
\section{Finite Multiple Sums and Mod $p$ Results}
\par
In this section we consider the finite sums
$$
\zt_n(i_1,\dots,i_k)=
\sum_{n\ge n_1>n_2>\dots>n_k\ge 1}\frac1{n_1^{i_1}\cdots n_k^{i_k}}
$$
and
$$
\zts_n(i_1,\dots,i_k)=
\sum_{n\ge n_1\ge n_2\ge\dots\ge n_k\ge 1}\frac1{n_1^{i_1}\cdots n_k^{i_k}} ;
$$
the multiple zeta values are
$$
\zt(i_1,\dots,i_k)=\lim_{n\to\infty}\zt_n(i_1,\dots,i_k),
$$
when the limit exists (i.e., when $i_1>1$).
\par
Let $\rho_n=e T \phi_n$, where
$\phi_n$ is the map defined in \S2, $T$ is the automorphism
of $\QSym$ sending $M_I$ to $M_{\bar I}$, and $e$ is the function
that sends $t_i$ to $\frac1{i}$.  Then $\rho_n:(\H^1,*)\to\mathbf R$ 
is a homomorphism sending $z_{i_1}z_{i_2}\cdots z_{i_k}$ in $\H^1$
to $\zt_n(i_1,\dots,i_k)$.
We can combine the homomorphisms $\rho_n$ into a homomorphism
$\rho$ that sends $w\in\H^1$ to the real-valued sequence $n\to\rho_n(w)$.
We shall write $\zt(I)$ for the real-valued sequence $n\to \zt_n(I)$
(and similarly for $\zts(I)$), so $\rho\phi^{-1}$ sends $M_I$ to $\zt(I)$ 
and $E_I$ to $\zts(I)$.  For example, we can apply $\rho\phi^{-1}$
to equation (\ref{sse}) above to get
\begin{equation}
\label{AS}
\zt(I)=\sum_{J\preceq I}(-1)^{\ell(I)-\ell(J)}\zts(J) .
\end{equation}
\par
Applying the homomorphism $\rho_n\phi^{-1}$ to Theorem 2.3 and 
equation (\ref{EM}), we get formulas for symmetric sums 
of $\zt_n(I)$ and $\zts_n(I)$ in terms of length one sums $\zt_n(m)$
(cf. Theorems 2.1 and 2.2 of \cite{H1}).
\begin{thm} For any composition $I=(i_1,\cdots,i_k)$,
\begin{align*}
\sum_{\si\in\Si_k}\zts_n(\si\cdot I)&=
\sum_{\text{partitions $\mathcal B=\{B_1,\dots,B_l\}$ of $\{1,\dots,k\}$}}
c(\mathcal B)\zt_n(b_1)\cdots \zt_n(b_l)\\
\sum_{\si\in\Si_k}\zt_n(\si\cdot I)&=
\sum_{\text{partitions $\mathcal B=\{B_1,\dots,B_l\}$ of $\{1,\dots,k\}$}}
(-1)^{k-l}c(\mathcal B)\zt_n(b_1)\cdots \zt_n(b_l) ,
\end{align*}
where $c(\mathcal B)=(\card B_1-1)!\cdots (\card B_l-1)!$ and
$b_s=\sum_{j\in B_s}i_j$.
\end{thm}
\par
We consider two operators on the space $\RN$ of real-valued sequences.
First, there is the partial-sum operator $\Si$, given by
$$
(\Si a)_n = \sum_{i=0}^n a_i 
$$
for $a\in\RN$.
Second, there is the operator $\del$ given by
$$
(\del a)_n = \sum_{i=0}^n \binom{n}{i}(-1)^i a_i .
$$
It is easy to show that $\Si$ and $\del$ generate a dihedral group
within the automorphisms of $\RN$, i.e., $\del^2=\id$ and $\Si\del=
\del\Si^{-1}$.  It follows that $(\Si\del)^2=\id$.  
We have the following result on multiple sums.
\begin{thm} For any composition $I$, $\Si\del\zts(I)=-\zts(I^*)$.
\end{thm}
\begin{proof} We proceed by induction on $|I|$.  The weight one
case is $\Si\del\zts(1)=\del\Si^{-1}\zts(1)=-\zts(1)$, i.e.,
$$
\sum_{k=1}^n\frac{(-1)^k}{k}\binom{n}{k}=-\sum_{k=1}^n\frac1{k} ,
$$
a classical (but often rediscovered) formula which actually goes 
back to Euler \cite{Eu}.  For $I=(i_1,i_2,\dots,i_k)$,
it is straightforward to show that $\del\zts_n(I)=\frac1{n}\del f(n)$, 
where $f\in\RN$ is given by
$$
f_n=\begin{cases}\zts_n(i_2,\dots,i_k),&\text{if $i_1=1$;}\\
\Si^{-1}\zts_n(i_1-1,i_2,\dots,i_k),&\text{otherwise.}\end{cases}
$$
\par
Now suppose the theorem has been proved for all $I$ of weight less
than $n$, and let $I=(i_1,\dots,i_k)$ have weight $n$.  There are
two cases:  $i_1=1$, and $i_1>1$.  In the first case, let
$(i_2,\dots,i_k)^*=J=(j_1,\dots,j_r)$.  By the assertion of the
preceding paragraph and the induction hypothesis,
$$
\Si\del\zts_n(I)=\Si\left(\frac1{n}\del\zts_n(J^*)\right)
=-\Si\left(\frac1{n}\Si^{-1}\zts_n(J)\right)
=-\zts_n(j_1+1,j_2,\dots,j_r).
$$
But evidently $I^*=(j_1+1,j_2,\dots,j_r)$, so the theorem holds in 
this case.
\par
If $i_1>1$, we instead write $(i_1-1,i_2,\dots,i_k)^*=J=(j_1,\dots,j_r)$.
Then
\begin{multline*}
\Si\del\zts_n(I)=\Si\left(\frac1{n}\del\Si^{-1}\zts_n(J^*)\right)
=\Si\left(\frac1{n}\Si\del\zts_n(J^*)\right)
=-\Si\left(\frac1{n}\zts_n(J)\right)\\
=-\zts_n(1,j_1,\dots,j_r).
\end{multline*}
But $I^*=(1,j_1,\dots,j_r)$, so the theorem holds in this case as well.
\end{proof}
\par
The proof of the preceding result is essentially a formalization of
the procedure in App. B of \cite{Ve}.  Recalling the automorphism
$\psi$ of $\H^1$ defined by equation (\ref{psi}), we note that
Theorem 4.2, together with part (2) of Theorem 3.2, says that the 
diagram
\begin{equation}
\begin{CD}
\QSym @>{\psi}>>\QSym\\
@V{\rho}VV @V{\rho}VV\\
\RN @>{\Si\del}>>\RN
\label{cmdiag}
\end{CD}
\end{equation}
commutes.
\par
We now turn to mod $p$ results about $\zts_{p-1}(I)$ and $\zt_{p-1}(I)$, 
where $p$ is a prime. 
Since the sums $\zt_{p-1}(I)$ and $\zts_{p-1}(I)$ contain no 
factors of $p$ in the denominators, they can be regarded as elements 
of the field $\Zp$.  
The following result about length one harmonic sums is
well known (cf. \cite{HW}, pp. 86-88).
\begin{thm} 
$\zt_{p-1}(k)\equiv 0 \mod p$ for all prime $p>k+1$.
\end{thm}
\begin{proof}
Let $S=\{1,2\dots,p-1\}\subset\Zp$ be the group of units.
For $c\in S$ we have
\[
c^k\zt_{p-1}(k)=\sum_{a\in S}c^ka^{-k}=\sum_{a\in S}(c^{-1}a)^{-k}=\zt_{p-1}(k)
\quad\text{in}\quad\Zp
\]
since multiplication by any unit of $\Zp$ permutes $S$, so 
\[
(c^k-1)\zt_{p-1}(k)=0\quad\text{in}\quad\Zp .
\]
Choosing $c\ne 1$, we see that $\zt_{p-1}(k)=0$ in $\Zp$.
\end{proof}
Combining the preceding result with Theorem 4.1 gives the following.
\begin{thm} For any composition $I=(i_1,\dots,i_k)$ and 
prime $p>|I|+1$,
$$
\sum_{\si\in\Si_k}\zt_{p-1}(\si\cdot I)\equiv
\sum_{\si\in\Si_k}\zts_{p-1}(\si\cdot I)\equiv 0\mod p .
$$
\end{thm}
In particular it follows that, for $I=(r,r,\dots,r)$ ($k$ repetitions), 
we have
\begin{equation}
\zt_{p-1}(I)\equiv \zts_{p-1}(I)\equiv 0\mod p
\label{vanish}
\end{equation}
for prime $p>rk+1$ (cf. Theorem 2.14 of \cite{Zh}).  There is
the following result relating sums associated to $I$ and $\bar I$
(cf. Lemma 3.3 of \cite{Zh}).
\begin{thm} For any composition $I$, $\zt_{p-1}(I)\equiv 
(-1)^{|I|}\zt_{p-1}(\bar I)\mod p$, and similarly
$\zts_{p-1}(I)\equiv(-1)^{|I|}\zts_{p-1}(\bar I)\mod p$.
\end{thm}
\begin{proof} Let $I=(i_1,\dots,i_k)$.  Working mod $p$, we have
\begin{multline*}
\zt_{p-1}(I)\equiv\sum_{p>a_1>\dots>a_k>0}\frac1{a_1^{i_1}\cdots a_k^{i_k}}
\equiv\sum_{p>a_1>\dots>a_k>0}\frac{(-1)^{i_1+\dots+i_k}}
{(p-a_1)^{i_1}\cdots (p-a_k)^{i_k}}\\
\equiv\sum_{0<b_1<\dots<b_k<p}\frac{(-1)^{i_1+\dots+i_k}}
{b_1^{i_1}\cdots b_k^{i_k}}=(-1)^{|I|}\zt_{p-1}(\bar I),
\end{multline*}
and similarly for $\zts_{p-1}$.
\end{proof}
An immediate consequence is that $\zts_{p-1}(I)\equiv \zt_{p-1}(I)
\equiv 0\mod p$ if $I=\bar I$ and $|I|$ is odd.  
Another consequence is that
$\zts_{p-1}(i,j)\equiv \zt_{p-1}(i,j)\equiv 0\mod p$ when $p>i+j+1$
and $i+j$ is even.  This is because
$$
\zts_{p-1}(i,j)+\zts_{p-1}(j,i)\equiv 0\mod p
$$
for $p>i+j+1$ by Theorem 4.1, while $\zts_{p-1}(i,j)\equiv \zts_{p-1}(j,i)
\mod p$ when $i+j$ is even by Theorem 4.5.
\par
We have the following result relating $\zts_{p-1}(I)$ and $\zts_{p-1}(I^*)$.
\begin{thm} $\zts_{p-1}(I)\equiv -\zts_{p-1}(I^*)\mod p$ for all primes $p$.
\end{thm}
\begin{proof} Let $f$ be a sequence.  From the definition of $\del$
$$
\Si\del f(n)=\sum_{i=0}^n\binom{n+1}{i+1}(-1)^if(i),
$$
so taking $n=p-1$ gives
$$
\Si\del f(p-1)\equiv (-1)^{p-1}f(p-1)\equiv f(p-1)\mod p .
$$
Now take $f=\zts(I)$ and apply Theorem 4.2.
\end{proof}
\par
Define, for each prime $p$, a map $\chi_p:\H^1\to\Zp$ by
$\chi_p(w)=\rho_{p-1}(w)$.  Then we can use the commutative diagram 
(\ref{cmdiag}) to restate Theorem 4.6 in the following ``algebraic'' form,
which corresponds to the duality theorem for multiple zeta values
as formulated in Corollary 6.2 of \cite{H2}.
\begin{thm} 
As elements of $\Zp$, $\chi_p(w)=\chi_p(\psi(w))$ for words $w$ of $\H^1$.
\end{thm}
For example, since $\psi(x^2y^3)=-x^2y^3-xy^4-yxy^3-y^5$, we have
$$
\zt_{p-1}(3,1,1)\equiv -\zt_{p-1}(3,1,1)-\zt_{p-1}(2,1,1,1)
-\zt_{p-1}(1,2,1,1)-\zt_{p-1}(1,1,1,1,1) \mod p .
$$
(In fact, for $p>6$ it follows from Theorem 7.1 below that both sides 
are congruent mod $p$ to $\frac12B_{p-5}$, where $B_i$ is the $i$th 
Bernoulli number.)
\section{Harmonic Sums of Height One}
Following the terminology of \cite{H4}, we say a word of $\H^1$ of the 
form $x^{h-1}y^k$ has height one.  We shall call the corresponding 
multiple harmonic sums $\zt_n(h,1,\dots,1)$ sums of height one.
Harmonic sums of this form have a number of special properties,
which we discuss in this section.
\par
In this section only we use superscripts for repetition in compositions, 
so $(n,1^k)$ means the composition of weight $n+k$ with $k$ repetitions of 1.
For compositions of this form, we have the following result relating the 
two kinds of harmonic sums mod $p$.
\begin{thm} For any prime $p>k$,
$$
\zts_{p-1}(h,1^{k-1})\equiv (-1)^h \zt_{p-1}(h,1^{k-1}) \mod p .
$$
\end{thm}
\begin{proof}
Equate parts (1) and (2) of Theorem 3.1 and apply the homomorphism 
$\rho_{p-1}$ to get
\begin{equation}
\label{ant12}
(-1)^{\ell(I)}\zts_{p-1}(\bar I)=\sum_{I_1\sqcup\dots\sqcup I_l=I}
(-1)^l\zt_{p-1}(I_1)\cdots \zt_{p-1}(I_l)
\end{equation}
for any $p$.  Now set $I=(1^{k-1},h)$ in equation (\ref{ant12}) and reduce 
both sides mod $p$.  The hypothesis insures that all terms on the 
right-hand side are zero mod $p$ except the one with $l=1$, and so we have
$$
(-1)^{k-1}\zts_{p-1}(h,1^{k-1})\equiv\zt_{p-1}(1^{k-1},h) \mod p ,
$$
and the conclusion follows by Theorem 4.5.
\end{proof}
Combining the preceding result with the duality theorem gives
the following congruence.  It can be compared with Theorem 4.4 of \cite{H1},
which asserts that $\zt(h+1,1^{k-1})=\zt(k+1,1^{h-1})$ for positive
integers $h$ and $k$.
\begin{thm} If $p$ is a prime with $p>\max\{k,h\}$, then
$$
\zt_{p-1}(h,1^{k-1})\equiv \zt_{p-1}(k,1^{h-1})\mod p .
$$
\end{thm}
\begin{proof} First note that $(h,1^{k-1})^*=(1^{h-1},k)$.  So, combining
Theorems 4.6 and 4.5,
\begin{equation}
\zts_{p-1}(h,1^{k-1})\equiv -\zts_{p-1}(1^{h-1},k)
\equiv (-1)^{h+k}\zts_{p-1}(k,1^{h-1})\mod p .
\label{S11}
\end{equation}
Now by the preceding result 
\begin{align*}
\zts_{p-1}(h,1^{k-1})&\equiv (-1)^h \zt_{p-1}(h,1^{k-1}) \mod p , \\
\zts_{p-1}(k,1^{h-1})&\equiv (-1)^k \zt_{p-1}(k,1^{h-1}) \mod p, 
\end{align*}
so the conclusion follows from congruence (\ref{S11}).
\end{proof}
\par
The height one harmonic sums can be written in terms of Stirling numbers.  
The (unsigned) Stirling numbers of the first kind (Stirling cycle numbers) 
are given by
$$
\stone{n}{k}= \text{number of permutations of $\{1,2,\dots,n\}$ with $k$ 
disjoint cycles,}
$$
and the Stirling number of the second kind (Stirling subset numbers)
are given by
$$
\sttwo{n}{k} = \text{number of partitions of $\{1,2,\dots,n\}$ with $k$ 
blocks.}
$$
For basic properties of Stirling numbers see, e.g., \cite{GKP}.
Note that
$$
\stone{n+1}{k+1}=e_{n-k}(1,2,\dots,n),
$$
since each side is the coefficient of $x^k$ in the expansion of $x(x+1)\cdots
(x+n)$.  Hence
\begin{equation}
\stone{n+1}{k+1}=n!\sum_{n\ge j_1>\dots>j_{n-k}\ge 1}
\frac{j_1\cdots j_{n-k}}{n!}
=n!\sum_{n\ge i_1>\dots>i_k\ge 1}\frac1{i_1\cdots i_k}=n!\zt_n(1^k) .
\label{st1}
\end{equation}
An immediate consequence is the following result, which expresses 
height one sums in terms of Stirling numbers of the first kind.
\begin{thm} For positive integers $h$, $k$, and $n$,
$$
\zt_n(h,1^{k-1})=\sum_{j=k}^n \frac{\stone{j}{k}}{j^{h-1}j!} .
$$
\end{thm}
\begin{proof} We have
$$
\zt_n(h,1^{k-1})=\sum_{n\ge i_1>\dots>i_k\ge 1}\frac1{i_1^hi_2\cdots i_k}
=\sum_{j=k}^n\frac1{j^h}\zt_{j-1}(1^{k-1}) =
\sum_{j=k}^n \frac{\stone{j}{k}}{j^h(j-1)!}
$$
(where we used equation (\ref{st1}) in the last step), from which the
conclusion follows.
\end{proof}
In doing computations, the following mod $p$ result is often 
useful.
\begin{thm} For prime $p$ and positive integers $h,k<p$,
$$
\zt_{p-1}(h,1^{k-1})\equiv \sum_{j=1}^{p-k} (-1)^j(-j)^{p-h}(j-1)!
\sttwo{p-k}{j} \mod p .
$$
\end{thm}
\begin{proof}
We start by setting $n=p-1$ and reversing the order of summation in 
the preceding result:
$$
\zt_{p-1}(h,1^{k-1})=\sum_{j=1}^{p-k}\frac{\stone{p-j}{k}}{(p-j)^{h-1}(p-j)!} .
$$
By Fermat's theorem, $(p-j)^{1-h}\equiv (p-j)^{p-h}\equiv (-j)^{p-j}$
mod $p$.  Also, since
$$
(-1)^{j-1}\equiv \binom{p-1}{j-1}\equiv\frac{(p-1)!}{(j-1)!(p-j)!} \mod p ,
$$
it follows that 
$$
\frac1{(p-j)!}\equiv (-1)^j(j-1)! \mod p
$$
since $(p-1)!\equiv -1$ mod $p$ by Wilson's theorem.  Hence
\begin{equation}
\zt_{p-1}(h,1^{k-1})\equiv 
\sum_{j=1}^{p-k}(-1)^j(j-1)!(-j)^{p-h}\stone{p-j}{k} \mod p .
\label{Ast2}
\end{equation}
The conclusion then follows from congruence (\ref{Ast2}) and
a mod $p$ relation between the two kinds of Stirling numbers:
$$
\stone{n}{k}\equiv \sttwo{p-k}{p-n} \mod p
$$
for prime $p$ and $1\le k\le n\le p-1$, which follows by induction
on $n$ using the recurrence relations for the two kinds of Stirling
numbers.
\end{proof}
\section{Results for Short Lengths}
In this section we prove some results about the sums $\zts_{p-1}(I)$
for $\ell(I)\le 5$, which will be useful for the calculations
of the next section. (Similar results hold for the $\zt_{p-1}(I)$,
but the $\zts_{p-1}(I)$ are more convenient in view of Theorem 4.6.)
Henceforth we will assume that all congruences are mod $p$.
The following result expresses
$\zts_{p-1}(I)$ in terms of a Bernoulli number when $\ell(I)=2$
(cf. \cite{Zh}, Theorem 3.1).
\begin{thm} For $i,j$ positive and prime $p>i+j+1$,
$$
\zts_{p-1}(i,j)\equiv\frac{(-1)^i}{i+j}\binom{i+j}{i}B_{p-i-j} .
$$
\end{thm}
\begin{proof}
We use the standard identity expressing sums of powers in terms
of Bernoulli numbers (see, e.g., \cite{IR}):
\begin{equation}
\sum_{a=1}^{n-1}a^r=\frac1{r+1}\sum_{k=0}^r\binom{r+1}{k}B_kn^{r+1-k} .
\label{bern}
\end{equation}
Using Fermat's theorem and equation (\ref{bern}), we have
\begin{multline*}
\zts_{p-1}(i,j)\equiv \zt_{p-1}(i,j)=\sum_{a=1}^{p-1}\frac1{a^i}
\sum_{b=1}^{a-1}\frac1{b^j}\equiv\sum_{a=1}^{p-1}\frac1{a^i}\sum_{b=1}^{a-1}
b^{p-1-j}\\
=\sum_{a=1}^{p-1}\frac1{a^j}\frac1{p-j}\sum_{k=0}^{p-1-j}\binom{p-j}{k}
B_ka^{p-j-k}=
\frac1{p-j}\sum_{k=0}^{p-1-j}\binom{p-j}{k}B_k\sum_{a=1}^{p-1}a^{p-i-j-k} .
\end{multline*}
Now $\sum_{a=1}^{p-1}a^{p-i-j-k}\equiv 0$ unless $k=p-i-j$.  So the
sum reduces to
$$
\zts_{p-1}(i,j)\equiv\frac1{p-j}\binom{p-j}{p-i-j}B_{p-i-j}(p-1)\equiv
\frac1{j}\binom{p-j}{i}B_{p-i-j},
$$
from which the conclusion follows.
\end{proof}
Since odd Bernoulli numbers are zero, this result implies our earlier
observation that $\zts_{p-1}(i,j)\equiv 0$ when $i+j$ is even.  Note
also that it implies $\zts_{p-1}(i,1)\equiv B_{p-1-i}$ for $i$ even,
so we can restate it as
\begin{equation}
\zts_{p-1}(i,j)\equiv \frac{(-1)^i}{i+j}\binom{i+j}{i}\zts_{p-1}(i+j-1,1),
\label{2to1}
\end{equation}
i.e., all double sums of odd weight $n$ can be written in terms 
of $\zts_{p-1}(n-1,1)$.  In fact, if the weight is odd  we can write 
all triple sums in terms of the same quantity (cf. \cite{Zh}, Theorem
3.5).
\begin{thm} If $n=i+j+k$ is odd and $p>n+1$, then
$$
\zts_{p-1}(i,j,k)\equiv\frac{1}{2n}\left[(-1)^i\binom{n}{i}+
(-1)^{i+j}\binom{n}{i+j}\right]\zts_{p-1}(n-1,1).
$$
\end{thm}
\begin{proof}
From Theorem 3.1 it follows that
\begin{equation}
\zt_m(I)=\sum_{I_1\sqcup I_2\sqcup\cdots\sqcup I_l=I}(-1)^{\ell(I)-l}
\zts_m(\bar I_1)\zts_m(\bar I_2)\cdots \zts_m(\bar I_l) 
\label{ASprod}
\end{equation}
for any composition $I$ and positive integer $m$.  In particular,
for $\ell(I)=3$ we have
$$
\zt_m(i,j,k)=\zts_m(i)\zts_m(j)\zts_m(k)-\zts_m(j,i)\zts_m(k)-
\zts_m(i)\zts_m(k,j)+\zts_m(k,j,i) .
$$
But also, from equation (\ref{AS}),
$$
\zt_m(i,j,k)=\zts_m(i,j,k)-\zts_m(i+j,k)-\zts_m(i,j+k)+\zts_m(i+j+k) .
$$
In the case $m=p-1$, this gives the congruence
$$
\zts_{p-1}(k,j,i)\equiv \zts_{p-1}(i,j,k)-\zts_{p-1}(i+j,k)-\zts_{p-1}(i,j+k) ,
$$
using Theorem 4.3, and if in addition $i+j+k$ is odd, we have
$$
2\zts_{p-1}(i,j,k)\equiv\zts_{p-1}(i+j,k)+\zts_{p-1}(i,j+k) 
$$
from Theorem 4.5.  Now use congruence (\ref{2to1}).
\end{proof}
For a composition
$I=(i_1,i_2,\dots,i_k)$ of odd weight and length $k>1$ it is convenient 
to define
\begin{equation}
C(I)=\sum_{j=1}^{k-1}(-1)^{i_1+\dots+i_j}\binom{|I|}{i_1+\dots+i_j} ,
\label{cI}
\end{equation}
so Theorem 6.2 says that
$$
\zts_{p-1}(I)\equiv \frac1{2n}C(I)\zts_{p-1}(n-1,1)
$$
when $n=|I|$ is odd and $\ell(I)=3$.  It is easy to show that
the function $C$ mirrors the properties of $\zts_{p-1}$ with
regard to reversal and duality, i.e., $C(\bar I)=-C(I)$ and
$C(I^*)=-C(I)$.
\par
In the even-weight case, we can express all length-4 sums in terms
of triple sums and products of lower-weight sums as follows.
\begin{thm} If $n=i+j+k+l$ is even and $p>n+1$, then
\begin{multline*}
2\zts_{p-1}(i,j,k,l)\equiv \zts_{p-1}(i+j,k,l)+\zts_{p-1}(i,j+k,l)\\
+\zts_{p-1}(i,j,k+l)+\zts_{p-1}(i,j)\zts_{p-1}(k,l) .
\end{multline*}
\end{thm}
\begin{proof}
From equation (\ref{ASprod}) with $\ell(I)=4$,
\begin{multline*}
\zt_m(i,j,k,l)=-\zts_m(l,k,j,i)+\zts_m(i)\zts_m(l,j,k)+
\zts_m(j,i)\zts_m(l,k)+\zts_m(k,j,i)\zts_m(k)\\
-\zts_m(i)\zts_m(j)\zts_m(l,k)
-\zts_m(i)\zts_m(k,j)\zts_m(l)
-\zts_m(j,i)\zts_m(k)\zts_m(l)
+\zts_m(i)\zts_m(j)\zts_m(k)\zts_m(l) .
\end{multline*}
Set $m=p-1$ and reduce mod $p$ to get
$$
\zt_{p-1}(i,j,k,l)\equiv -\zts_{p-1}(l,k,j,i)+\zts_{p-1}(j,i)\zts_{p-1}(l,k) .
$$
But equation (\ref{AS}) gives, after reduction mod $p$,
$$
\zt_{p-1}(i,j,k,l)\equiv \zts_{p-1}(i,j,k,l)-\zts_{p-1}(i+j,k,l)
-\zts_{p-1}(i,j+k,l)-\zts_{p-1}(i,j,k+l) ,
$$
and the conclusion follows by the use of Theorem 4.5.
\end{proof}
\par
If we specialize the preceding result to the case $k=i,l=j$, we
have
$$
2\zts_{p-1}(i,j,i,j)\equiv \zts_{p-1}(i+j,i,j)+\zts_{p-1}(i,i+j,j)+
\zts_{p-1}(i,j,i+j)+\zts_{p-1}(i,j)^2 .
$$
But the sum of the first three terms on the right-hand side is 
zero mod $p$ by Theorem 4.4, so 
\begin{equation}
2\zts_{p-1}(i,j,i,j)\equiv \zts_{p-1}(i,j)^2 .
\label{square}
\end{equation}
We note that we also have $\zt_{p-1}(i,j,i,j)\equiv \zts_{p-1}(i,j,i,j)$,
which follows from equation (\ref{AS}), Theorem 4.4, and the
comments following Theorem 4.5.
If $i$ and $j$ have the same parity, then $i+j$ is even and the right-hand
side of (\ref{square}) is zero mod $p$ by Theorem 4.5, so in this case
$$
\zts_{p-1}(i,j,i,j)\equiv \zt_{p-1}(i,j,i,j) \equiv 0 .
$$
Cf. Theorem 3.18 of \cite{Zh}.
\par
It is evident that the ideas of the two preceding results can be
used to write $\zts_{p-1}(I)$, where $|I|$ and $\ell(I)$ are of the same
parity, in terms of sums of shorter lengths, together with products
of sums of lower weight.  We state one more result of this type,
since we need it for our computations in weight 9; we omit the
proof, which is straightforward.  
\begin{thm} If $I=(i,j,k,l,m)$ is a composition of odd weight $n$ and 
$p>n+1$, then
\begin{multline*}
2\zts_{p-1}(I)\equiv \zts_{p-1}(i+j,k,l,m)+\zts_{p-1}(i,j+k,l,m)+
\zts_{p-1}(i,j,k+l,m)+\zts_{p-1}(i,j,k,l+m)\\
-\frac1{2n}C(I)\zts_{p-1}(n-1,1)+\zts_{p-1}(i,j)\zts_{p-1}(k,l,m)
+\zts_{p-1}(i,j,k)\zts_{p-1}(l,m) ,
\end{multline*}
where $C(I)$ is given by equation (\ref{cI}).
\end{thm}
\section{Calculations in Low Weights}
In this section, we find, for $n\le 9$, sets of quantities
that generate all multiple harmonic sums $\zts_{p-1}(I)$ 
of weight $n$ mod $p$ when $p>n+1$.
In fact, all our generators will be sums of the form
$\zts_{p-1}(2h,1,\dots,1)$ ($\equiv \zt_{p-1}(2h,1,\dots,1)$ by Theorem 5.1), 
or products of such sums.
To make the notation less cumbersome we will use an abbreviated 
form in referring to specific sums, e.g., we write $\zts_{221}$ instead 
of $\zts_{p-1}(2,2,1)$.  We will also assume that $p>|I|+1$ in
considering any sum $\zts_{p-1}(I)$.
\par
The congruence (\ref{vanish}) implies $\zts_{p-1}(I)\equiv 0$ when $I$ has 
length 1 or $|I|$, so the only sums $\zts_{p-1}(I)$ with $|I|\le 3$
that are (possibly) nonzero mod $p$ are $\zts_{21}\equiv B_{p-3}$ 
and $\zts_{12}\equiv -B_{p-3}$.
Now $p$ divides the numerator of $B_{p-3}$ only for certain irregular 
primes $p$, the only examples with $p<1.2\times 10^7$ being 
$p=16,843$ and $p=2,124,679$ \cite{J,BCEM,B}.
\par
For $I$ of weight 4 we have $\zts_{p-1}(I)\equiv 0$; this follows
from congruence (\ref{vanish}) if $\ell(I)=1$ or 4, by the remarks
following Theorem 4.5 if $\ell(I)=2$, and by duality if
$\ell(I)=3$ (since in that case $\ell(I^*)=2$).
\par
In weight 5, it suffices by duality to consider $\zts_{p-1}(I)$ with
$\ell(I)\le 3$.  Theorems 6.1 and 6.2 give the following result, which
implies that all weight 5 sums vanish mod $p$ if and only if $p$
divides the numerator of $B_{p-5}$; the only known $p$ for which this 
happens is $p=37$, and this is the only such $p<1.2\times 10^7$ \cite{B}.
\begin{thm}
All $\zts_{p-1}(I)$ with $|I|=5$ are multiples of $\zts_{41}\equiv B_{p-5}$.  
In particular, $\zts_{32}\equiv -2\zts_{41}$, $\zts_{311}\equiv-\frac12\zts_{41}$, 
and $\zts_{221}\equiv\frac32\zts_{41}$.
\end{thm}
For $|I|=6$, it follows from the remarks after Theorem 4.5 and duality
that $\zts_{p-1}(I)\equiv 0$ unless $\ell(I)=3$ or 4.  By duality it suffices
to consider length 3.  Here we have the following result, which implies
that all harmonic sums of weight 6 vanish if and only if $p$ divides
$B_{p-3}$.
\begin{thm} All the $\zts_{p-1}(I)$ with $|I|=6$ are multiples of
$\zts_{411}\equiv -\frac16B_{p-3}^2$.  In particular,
$\zts_{141}\equiv -2\zts_{411}$,
$\zts_{312}\equiv -\zts_{411}$,
$\zts_{321}\equiv -2\zts_{411}$, and
$\zts_{231}\equiv 3\zts_{411}$.
\end{thm}
\begin{proof} From Theorem 4.4, the symmetric sum 
$$
\zts_{411}+\zts_{141}+\zts_{114}\equiv2\zts_{411}+\zts_{141}
$$
is congruent to zero, from which the statement about $\zts_{141}$
follows.  For the same reason
$$
\zts_{3111}+\zts_{1311}+\zts_{1131}+\zts_{1113}\equiv 2\zts_{3111}+2\zts_{1311}
$$
is congruent to zero; applying duality, we get $\zts_{411}+\zts_{312}\equiv 0$.  
Theorem 4.4 also gives
\begin{equation}
\zts_{321}+\zts_{231}\equiv -\zts_{312}\equiv \zts_{411}
\label{wt6e1}
\end{equation}
Now multiplying in $\QSym$ and reducing mod $p$, we have
$$
0\equiv \zts_{1}\zts_{221}
\equiv\zts_{1221}+\zts_{2121}+2\zts_{2211}-\zts_{321}-\zts_{231}-\zts_{222}
\equiv \zts_{2121}+2\zts_{2211}-\zts_{321}-\zts_{231} .
$$
Apply duality (and Theorem 4.5) to obtain from this
\begin{equation}
3\zts_{321}+2\zts_{231}\equiv 0 ,
\label{wt6e2}
\end{equation}
and solve congruences (\ref{wt6e1}) and (\ref{wt6e2}) to obtain $\zts_{321}$ 
and $\zts_{231}$ in terms of $\zts_{411}$.
\par
Finally, $\zts_{2121}\equiv \frac12(\zts_{21})^2$ by congruence (\ref{square}),
while from duality we have $\zts_{2121}\equiv -\zts_{231}$.  Since we've
already shown $\zts_{231}\equiv 3\zts_{411}$, the congruence
$(\zts_{21})^2\equiv -6\zts_{411}$ follows.
\end{proof}
In weight 7, we need only consider sums of length 4 or less.
By Theorems 6.1 and 6.2 we can express all sums of
length 2 and 3 in terms of $\zts_{61}\equiv B_{p-7}$, as follows.
\begin{align}
\zts_{52}\equiv& -3\zts_{61}\\
\zts_{43}\equiv&\ 5\zts_{61}\\
\zts_{511}\equiv& -\zts_{61}\\
\zts_{421}\equiv&\ 3\zts_{61}\\
\zts_{412}\equiv&\ \zts_{61}\\
\zts_{241}\equiv&\ 2\zts_{61}\\
\zts_{331}\equiv& -2\zts_{61}\\
\zts_{322}\equiv& -4\zts_{61}
\end{align}
For the length 4 sums we have the following result.
\begin{thm} The length 4 harmonic sums of weight 7 can be written 
in terms of $\zts_{61}\equiv B_{p-7}$ and $\zts_{4111}$.  In particular
$\zts_{3121}\equiv \zts_{2221}\equiv \zts_{4111}$, and
\begin{align*}
\zts_{1411}&\equiv 3\zts_{61}-3\zts_{4111}\\
\zts_{1321}&\equiv 9\zts_{61}-3\zts_{4111}\\
\zts_{1312}&\equiv -9\zts_{61}+5\zts_{4111}\\
\zts_{3211}&\equiv -6\zts_{61}+\zts_{4111} .
\end{align*}
\end{thm}
\begin{proof} 
Duality gives $\zts_{1411}\equiv \zts_{3112}$, $\zts_{3121}\equiv \zts_{2311}$,
and $\zts_{2212}\equiv \zts_{1321}$.  Thus, it suffices to consider the
seven length-four sums $\zts_{4111}$, $\zts_{1411}$, $\zts_{3121}$,
$\zts_{2212}$, $\zts_{1312}$, $\zts_{2221}$, and $\zts_{3211}$.  We can
obtain relations among these sums and $\zts_{61}$ by multiplying lower-weight
sums in $\QSym$, reducing mod $p$, and then applying congruences (27-34)
above to get all terms of length three or less in terms of $\zts_{61}$.
For example, multiplying $\zts_1$ by $\zts_{411}$ gives
\begin{equation}
3\zts_{61}\equiv \zts_{1411}+3\zts_{4111} .
\end{equation}
Similarly, consideration of the products $\zts_2\zts_{311}$, 
$\zts_2\zts_{221}$, $\zts_{21}\zts_{31}$, and $\zts_{11}\zts_{32}$ gives
respectively
\begin{align}
-3\zts_{61}&\equiv 2\zts_{3121}+\zts_{3211}+\zts_{1411}\\
9\zts_{61}&\equiv 3\zts_{2221}+\zts_{2212}\\
-3\zts_{61}&\equiv -\zts_{1312}+3\zts_{3121}+2\zts_{3211}\\
-3\zts_{61}&\equiv \zts_{1411}+\zts_{2212}+\zts_{1312}+\zts_{3211} .
\end{align}
Finally, multiplying $\zts_1$ by $\zts_{3111}$ and using duality gives
\begin{equation}
-3\zts_{61}\equiv \zts_{4111}+\zts_{1411}+\zts_{3211}+\zts_{3121}
\end{equation}
The six congruences (35-40) can then be solved for $\zts_{61}$ and
$\zts_{4111}$ to obtain the conclusion.
\end{proof}
\begin{remark} 
In fact, there is an additional congruence
$\zts_{4111}\equiv \frac{27}{16}\zts_{61}$
found by Kh. Hessami Pilehrood, T. Hessami Pilehrood, and R.
Tauraso \cite{HPT}.  Hence all quantities in weight 7 are multiples of 
$\zts_{61}\equiv B_{p-7}$.
\end{remark}
\par
In weight 8, it suffices by duality to consider $I$ with $\ell(I)\le 4$.
By Theorem 6.3, we can write all length 4 sums in terms of
triple sums and $\zts_{41}\zts_{21}$, and of course all double sums are
zero.  So it is enough to consider length 3, where we have the following
result.
\begin{thm}
All length 3, weight 8 sums can be written in terms of $\zts_{611}$
and $\zts_{41}\zts_{21}\equiv B_{p-5}B_{p-3}$.  In particular,
\begin{align}
\zts_{521}&\equiv \zts_{611}+\zts_{41}\zts_{21}\\
\zts_{512}&\equiv \frac12(-7 \zts_{611}- \zts_{41}\zts_{21})\\
\zts_{431}&\equiv \frac12(-25\zts_{611}-9 \zts_{41}\zts_{21})\\
\zts_{413}&\equiv \frac12(5\zts_{611}+\zts_{41}\zts_{21})\\
\zts_{332}\equiv \zts_{422}&\equiv 10\zts_{611}+2\zts_{41}\zts_{21} .
\end{align}
\end{thm}
\begin{proof}
Consideration of the products $\zts_1\zts_{211111}\equiv 0$ and
$\zts_1\zts_{121111}\equiv 0$, using duality, gives respectively
\begin{align}
\zts_{611}+\zts_{521}+\zts_{431}+\zts_{341}+\zts_{251}+\zts_{161}&\equiv 0,\\
\zts_{521}+\zts_{512}+\zts_{422}+\zts_{332}+\zts_{242}+\zts_{152}&\equiv 0 .
\end{align}
Now from symmetry (Theorem 4.4), $\zts_{161}=
-2\zts_{611}$, $\zts_{251}\equiv -\zts_{521}-\zts_{512}$,
and $\zts_{341}\equiv -\zts_{431}-\zts_{413}$.  Using these facts, congruences
(46-47) reduce to 
\begin{align}
\zts_{413}&\equiv \zts_{611}-\zts_{512},\\
\zts_{332}&\equiv \zts_{422}.
\end{align}
Now consider the product $\zts_1\zts_{31111}\equiv 0$.  Using duality,
this gives
$$
5\zts_{611}+\zts_{512}\equiv \zts_{4211}+\zts_{3311}+\zts_{2311} .
$$
Using Theorem 6.3, the right-hand side can be expressed in
terms of triple sums and $\zts_{41}\zts_{21}$, and then applying the preceding
relations we obtain
\begin{equation}
\zts_{521}\equiv -6\zts_{611}-2\zts_{512}.  
\end{equation}
Applying similar techniques to the product $\zts_2\zts_{2111}$ and using 
congruences (48-50) gives congruences (41), (42), and (44).
\par
Finally, using duality on the product $\zts_2\zts_{2112}\equiv 0$ gives
$$
0\equiv 4\zts_{1421}+2\zts_{4112}+2\zts_{2312}+\zts_{1331} .
$$
Now rewrite this using Theorem 6.3 and apply the previously established 
relations to obtain (43); congruence (45) then follows.
\end{proof}
\begin{remark} It cannot be the case that $\zts_{611}=u\zts_{41}\zts_{21}$, where
$u$ is a unit of $\Zp$ for all $p>9$, since when $p=37$ we have 
$\zts_{611}\equiv 7$ and $\zts_{41}\equiv 0$.
\end{remark}
\par
For sums of weight 9, duality allows us to restrict attention to $\zts_{p-1}(I)$
with $\ell(I)\le 5$.  Using Theorem 6.4 (and Theorem 7.2), we can 
rewrite the length 5 sums in terms of length 4 sums, $\zts_{81}$, and
the product $\zts_{411}\zts_{21}$.  So it is enough to consider the sums 
$\zts_{p-1}(I)$ with $\ell(I)=4$, for which we have the following result.
\begin{thm} All length 4, weight 9 harmonic sums can be written in 
terms of the three 
quantities $\b$, $\c\equiv B_{p-9}$, and $\e\equiv -\frac16B_{p-3}^3$.  
In particular,
\begin{xalignat*}{2}
\zts_{1611}&\equiv-3\b+\frac{19}3\c &
\zts_{5211}&\equiv\b-\frac{310}{27}\c+\frac13\e\\
\zts_{5112}&\equiv-4\b+\frac{190}{27}\c-\frac13\e &
\zts_{5121}&\equiv\b-\frac{70}{27}\c+\frac13\e\\
\zts_{2511}&\equiv2\b-\frac{89}{27}\c-\frac13\e &
\zts_{2151}&\equiv-7\b+\frac{544}{27}\c-\frac13\e \\
\zts_{1521}&\equiv-3\b+\frac{137}9\c-\e &
\zts_{4311}&\equiv\b+\frac{13}3\c \\
\zts_{4131}&\equiv-6\b+\frac{56}3\c-\e &
\zts_{4113}&\equiv7\b-\frac{133}9\c+\e\\
\zts_{3411}&\equiv-2\b-\frac{20}9\c &
\zts_{3141}&\equiv8\b-\frac{182}9\c+\e\\
\zts_{1431}&\equiv4\b-\frac{64}3\c+\e &
\zts_{4221}&\equiv\b+\frac{82}{27}\c-\frac13\e\\
\zts_{4212}&\equiv-4\b+\frac{617}{27}\c-\frac23\e &
\zts_{4122}&\equiv3\b-\frac{370}{27}\c+\frac13\e\\
\zts_{2421}&\equiv2\b-\frac{124}{27}\c+\frac43\e &
\zts_{2412}&\equiv-2\b+\frac{110}9\c\\
\zts_{2241}&\equiv\frac{122}{27}\c-\frac23\e &
\zts_{3321}&\equiv-2\b+\frac{290}{27}\c-\frac23\e\\
\zts_{3312}&\equiv2\b-\frac{289}{27}\c+\frac13\e &
\zts_{3231}&\equiv-\frac{292}{27}\c+\frac43\e\\
\zts_{2331}&\equiv\frac19\c-\e &
\zts_{3213}&\equiv4\b-\frac{64}3\c\\
\zts_{3132}&\equiv-\frac{43}{27}\c-\frac23\e &
\zts_{2322}&\equiv-\frac{28}3\c\\
\zts_{3222}&\equiv-6\c . &&
\end{xalignat*}
\end{thm}
\begin{proof} There are 56 length-4 sums, but by reversal
(Theorem 4.5) we can reduce to 28 of them.  Next we go from 28 length-4 
sums to 14 length-4 sums, together with $\zts_{81}$, by considering various 
products in $\QSym$.  
For example, from the products $\zts_1\zts_{422}$, $\zts_2\zts_{421}$, 
and $\zts_2\zts_{412}$ respectively, we obtain the relations
\begin{align*}
\frac{23}{3}\zts_{81} &\equiv -\zts_{2241}+\zts_{4122}+\zts_{4212}+\zts_{4221} \\
\frac{73}{3}\zts_{81} &\equiv \zts_{2421}+2\zts_{4221}+\zts_{4212} \\
\frac{23}{3}\zts_{81} &\equiv \zts_{2412}+\zts_{4212}+2\zts_{4122}
\end{align*}
which we can use to write $\zts_{2241}$, $\zts_{2421}$, and $\zts_{2412}$
in terms of $\zts_{4221}$, $\zts_{4212}$, $\zts_{4122}$, and $\zts_{81}$.
Similarly, we eliminate $\zts_{1611}$, $\zts_{2511}$, $\zts_{2151}$,
$\zts_{1521}$, $\zts_{3411}$, $\zts_{3141}$, $\zts_{1431}$, $\zts_{2331}$,
$\zts_{3213}$, $\zts_{3132}$, and $\zts_{2322}$.  Now we consider the
eight products $\zts_1\zts_{311111}$, $\zts_1\zts_{221111}$, $\zts_1\zts_{131111}$,
$\zts_1\zts_{113111}$, $\zts_1\zts_{212111}$, $\zts_1\zts_{211211}$, $\zts_1\zts_{211121}$,
and $\zts_1\zts_{122111}$.  Using duality and rewriting in terms of our
14 length-4 sums, we have the relations
\begin{align*}
\frac{28}{3}\c&\equiv 2\zts_{6111}+\zts_{5121}+\zts_{5112}+\zts_{4131}+\zts_{4113}\\
\frac{11}{3}\c&\equiv \zts_{5211}+\zts_{4221}+\zts_{4212}-\zts_{3321}\\
\frac{29}{3}\c&\equiv -\zts_{5112}-2\zts_{4122}+\zts_{3312}\\
-\frac{61}{3}\c&\equiv -\zts_{4113}+\zts_{4122}+2\zts_{3312}\\
\frac{4}{3}\c&\equiv \zts_{4311}-\zts_{4221}+\zts_{3321}+\zts_{3312}\\
11\c&\equiv -\zts_{4113}+\zts_{3141}+\zts_{4311}+\zts_{3321}+\zts_{3231}+\zts_{4122}+\zts_{4212}
+\zts_{4221}\\
-13\c&\equiv \zts_{5211}+\zts_{5112}+\zts_{5121}+\zts_{4122}-\zts_{4221}-\zts_{3231}-\zts_{3321}
-\zts_{3312}\\
-6\c&\equiv \zts_{3222} .
\end{align*}
We can use these congruences to solve for all the length-4 sums in terms
of the seven quantities $\c$, $\b$, $\zts_{5211}$, $\zts_{5112}$, $\zts_{5121}$,
$\zts_{4311}$, and $\zts_{4221}$.  
\par
From duality, $\zts_{41211}\equiv \zts_{33111}$.  Expanding out both
sides using Theorem 6.4, we obtain the congruence
$$
\zts_{5211}+\zts_{4311}+\zts_{4131}+\zts_{4122}\equiv
\zts_{6111}+\zts_{3411}+\zts_{3321}+\zts_{3312},
$$
which (after rewriting in terms of our seven chosen quantities),
allows us to solve for 
$\zts_{4221}$ in terms of the remaining six quantities.  Now apply
Theorem 6.4 similarly to the duality relations
$\zts_{41112}\equiv \zts_{15111}$, $\zts_{12411}\equiv \zts_{31122}$, and 
$\zts_{21411}\equiv \zts_{31131}$ to obtain respectively
\begin{align}
-\zts_{4311}+2\zts_{5121}+\zts_{5211}-2\zts_{61111}+21\zts_{81}&\equiv \zts_{21}\zts_{411}\\
3\zts_{4311}-\zts_{5121}+4\zts_{5211}-6\zts_{61111}+\frac{91}{3}\zts_{81}&\equiv \zts_{21}\zts_{411}\\
-3\zts_{4311}+2\zts_{5121}+\zts_{5211}+\frac{89}{3}\zts_{81}&\equiv \zts_{21}\zts_{411} .
\end{align}
Add 3 times congruence (51) to congruence (52), and compare the
resulting expression for $\zts_{21}\zts_{411}$ with that obtained by
adding congruence (53) to (52); the result is
\begin{equation}
\zts_{5121}\equiv \zts_{5211}+\frac{80}9\zts_{81}.
\label{sim1}
\end{equation}
Next substitute congruence (\ref{sim1}) into the result of subtracting
congruence (52) from (51) to get
\begin{equation}
\zts_{4311}\equiv \zts_{6111}+\frac{13}{3}\zts_{81} .
\label{sim2}
\end{equation}
Using congruences (\ref{sim1}) and (\ref{sim2}), any of 
the congruences (51-53) can be used to write $\zts_{5211}$ in
terms of $\zts_{6111}$, $\zts_{81}$, and $\zts_{21}\zts_{411}$.
\par
Finally, consideration of the product
$\zts_{11}\zts_{2221}$ gives the relation
\begin{multline*}
3\zts_{3231}+3\zts_{3321}+3\zts_{4221}+2\zts_{3222}+2\zts_{2331}+\zts_{2241}+\zts_{2421}\\
\equiv 3\zts_{31221}+3\zts_{32121}+2\zts_{32211}+2\zts_{22131}+\zts_{21231}+\zts_{23121} ,
\end{multline*}
which can be used to write $\zts_{5112}$ in terms of $\zts_{6111}$, $\zts_{81}$,
and $\zts_{21}\zts_{411}$.
\end{proof}
\begin{remarks} 1. The only primes $p<1.2\times 10^7$ with $B_{p-9}\equiv 0$
are $p=67$ and $p=877$ \cite{B}.  Here is a table of the values of our
generators at these primes, and at $p=16,843$, computed using \it Mathematica
\rm and Theorem 5.4.
\par
\begin{center}
\begin{tabular}{|l|c|c|c|}\hline
$p$ & 67 & 877 & 16,843 \\ \hline
$\zts_{6111}$ & 7 & 253 & 16,690 \\ \hline
$\zts_{81}$ & 0 & 0 & 14,820 \\ \hline
$\zts_{21}\zts_{411}$ & 4 & 354 & 0 \\ \hline
\end{tabular}
\end{center}
2. J. Zhao conjectured that
$$
\zts_{6111}\equiv -\frac19 \zts_{21}\zts_{411} + \frac{1889}{648} \zts_{81}
$$
(see \cite[Prop. 3.12]{Zh}).  This conjecture is consistent with the
data in the preceding remark.
\par\noindent
3. The congruence of the preceding remark is indeed true
\cite{HPT}.
\end{remarks}


\begin{thebibliography}{99}
\bibitem{BK}
J. Bl\"umlein and S. Kurth, Harmonic sums and Mellin transforms
up to two-loop order, \emph{Phys. Rev. D} {\bf 60} (1999), art. 01418.
\bibitem{BB}
D. Bowman and D. M. Bradley, Multiple polylogarithms:  a brief survey,
in \emph{$q$-Series with Applications to Combinatorics, Number Theory,
and Physics}, Contemp. Math., Vol. 291, American Mathematical Society, 
Providence, 2001, pp. 71-92.
\bibitem{B}
J. Buhler, personal communication.
\bibitem{BCEM}
J. Buhler, R. Crandall, R. Ernvall, and T. Mets\"ankyl\"a,
Irregular primes and cyclotomic invariants up to four million,
\emph{Math. Comp.} {\bf 61} (1993), 151-153.
\bibitem{Eh}
R. Ehrenborg, On posets and Hopf algebras, \emph{Adv. Math.} {\bf 119}
(1996), 1-25.
\bibitem{Eu}
L. Euler, Demonstratio insignis theorematis numerici circa unicias potestatum
binomialium, \emph{Nova Acta Acad. Sci. Petropol.} {\bf 15} (1799/1802), 33-43;
reprinted in \emph{Opera Omnia}, Ser. I, Vol. 16(2), B. G. Teubner, Leipzig,
1935, pp. 104-116.
\bibitem{Gei}
L. Geissinger, Hopf algebras of symmetric functions and class functions,
in \emph{Combinatoire et repr\'esentation de groupe symm\'etrique
(Strasbourg, 1976)}, Lecture Notes in Math., Vol. 579, Springer-Verlag,
New York, 1977, pp. 168-181.
\bibitem{Gel}
I. M. Gelfand, D. Krob, A. Lascoux, B. Leclerc, V. S. Retakh, and
J.-Y. Thibon, Noncommutative symmetric functions, \emph{Adv. Math.}
{\bf 112} (1995), 218-348.
\bibitem{Ges}
I. M. Gessel, ``Multipartite P-partitions and inner products of skew
Schur functions,'' in \emph{Combinatorics and Algebra}, Contemp. Math.,
Vol. 34, American Mathematical Society, Providence, 1984, pp. 289-301.
\bibitem{G}
A. G. Goncharov, Multiple polylogarithms and mixed Tate motives, 
{\tt arXiv math/010359}.
\bibitem{GKP}
R. L. Graham, D. E. Knuth, and O. Patashnik, \emph{Concrete Mathematics},
2nd ed., Addison-Wesley, Reading, Mass., 1994.
\bibitem{HW}
G. H. Hardy and E. M. Wright, \emph{An Introduction to the Theory
of Numbers}, 4th ed., Oxford University Press, London, 1960.
\bibitem{Hz}
M. Hazewinkel, The algebra of quasi-symmetric functions is free over
the integers, \emph{Adv. Math.} {\bf 164} (2001), 283-300.
\bibitem{HPT}
Kh. Hessami Pilehrood, T. Hessami Pilehrood, and R. Tauraso, New 
properties of multiple harmoinc sums modulo $p$ and $p$-analogues 
of Leschinger's series, \emph{Trans. Amer. Math. Soc.} {\bf 366} (2014),
3131-3159.
\bibitem{H1}
M. E. Hoffman, Multiple harmonic series, \emph{Pacific J. Math.} {\bf 152}
(1992), 275-290.
\bibitem{H2}
M. E. Hoffman, The algebra of multiple harmonic series, \emph{J. Algebra}
{\bf 194} (1997), 477-495.
\bibitem{H3}
M. E. Hoffman, Quasi-shuffle products, \emph{J. Algebraic Combin.}
{\bf 11} (2000), 49-68.
\bibitem{H4}
M. E. Hoffman, Algebraic aspects of multiple zeta values, in
\emph{Zeta Functions, Topology and Quantum Physics,} T. Aoki {\it et. al.} 
(eds.), Springer, New York, 2005, pp. 51-74.
\bibitem{HO}
M. E. Hoffman and Y. Ohno, Relations of multiple zeta values and their 
algebraic expression, \emph{J. Algebra} {\bf 262} (2003), 332-347.
\bibitem{J}
W. Johnson, Irregular primes and cyclotomic invariants, \emph{Math. Comp.} 
{\bf 29} (1975), 113-120.
\bibitem{IR}
K. Ireland and M. Rosen, \emph{A Classical Introduction to Modern Number
Theory}, Spring-Verlag, New York, 1982.
\bibitem{Ka}
C. Kassel, \emph{Quantum Groups}, Springer-Verlag, New York, 1995.
\bibitem{Mac} 
I. G. Macdonald, \emph{Symmetric Functions and Hall Polynomials,}
2nd ed., Oxford University Press, New York, 1995.
\bibitem{MR}
C. Malvenuto and C. Reutenauer, Duality between quasi-symmetric functions
and the Solomon descent algebra, \emph{J. Algebra} {\bf 177} (1995), 967-982.
\bibitem{MM}
J. W. Milnor and J. C. Moore, On the structure of Hopf algebras,
\emph{Ann. of Math. (2)} {\bf 81} (1965), 211-261.
\bibitem{St}
R. Stanley, \emph{Enumerative Combinatorics}, Vol. 1, Wordsworth and
Brooks/Cole, Monterey, California, 1986.
\bibitem{S}
M. Sweedler, \emph{Hopf Algebras}, W. A. Benjamin, Inc., New York, 1969.
\bibitem{T}
T. Terasoma, Mixed Tate Motives and multiple zeta values, {\emph Inv.
Math.} {\bf 199} (2002), 339-369.
\bibitem{Ve}
J. A. M. Vermaseren, Harmonic sums, Mellin transforms and integrals,
\emph{Int. J. Mod. Phys. A} {\bf 14} (1999), 2037-2076.
\bibitem{W}
M. Waldschmidt, Valeurs z\^eta multiples.  Une introduction,
\emph{J. Th\'eor. Nombres Bordeaux} {\bf 12} (2000), 581-595.
\bibitem{Z}
D. Zagier, Values of zeta functions and their applications, in
\emph{First European Congress of Mathematics, Vol. II (Paris, 1992)}, 
Birkh\"auser, Boston, 1994, pp. 497-512.
\bibitem{Zh} J. Zhao, Wolstenholme type theorem for multiple harmonic
sums, \emph{Int. J. Number Theory} {\bf 4} (2008), 73-106; preprint
\tt math.NT/0301252\rm .
\end{thebibliography}
\end{document}